\documentclass[11pt,a4paper,leqno,twoside, headinclude]{amsart}

\usepackage[tracking=true]{microtype}
\DeclareMicrotypeSet*[tracking]{my}%
  { font = */*/*/sc/* }%
\SetTracking{ encoding = *, shape = sc }{ 45 }

\title{Deciding almost freeness of an action is NP-hard}
\def\titl{Deciding almost freeness of an action is NP-hard}
\def\auth{Manuel Amann}
\date{August 16th, 2015}

\thanks{}

\author{\auth}

\usepackage[english]{babel}

\usepackage[colorlinks,pdftex, plainpages=false]{hyperref}
\hypersetup{pdftitle=\titl, pdfauthor=\auth, pdftoolbar=false,
plainpages=false, hyperindex=true, pdfdisplaydoctitle=true}

\hypersetup{colorlinks,%
citecolor=black,%
filecolor=black,%
linkcolor=black,%
urlcolor=black,%
pdftex}

\usepackage{color}

\usepackage{amsmath, amssymb, amscd, amsthm}
\usepackage{stmaryrd}
\usepackage{latexsym}
\usepackage[all]{xy}
\usepackage{pb-diagram}
\usepackage{rotating}
\usepackage{multicol}
\usepackage{lscape}
\usepackage{wasysym}
\usepackage{longtable}
\usepackage{enumerate}

\xyoption{all}

\newtheorem{theo}{Theorem}[section]
\newtheorem{main}{Theorem}
\newtheorem{maincor}[main]{Corollary}
\newtheorem*{main*}{Theorem}

\newtheorem*{mainprop*}{Proposition}

\newtheorem{mainconj}{Conjecture}

\newtheorem{defi2}[theo]{Definition}
\newtheorem*{defi2*}{Definition}

\newenvironment{defi*}{\begin{defi2*}\normalfont}{\end{defi2*}}

\newenvironment{defin*}[1]{\begin{defi2*}[#1]\normalfont}{\end{defi2*}}
\newtheorem*{rem2*}{Remark}
\newenvironment{rem*}{\begin{rem2*}\normalfont}{\hfill$\boxbox$\end{rem2*}}
\newtheorem{rem2}[theo]{Remark}
\newenvironment{rem}{\begin{rem2}\normalfont}{\hfill$\boxbox$\end{rem2}}

\newtheorem{lemma}[theo]{Lemma}

\newtheorem*{cor*}{Corollary}

\newtheorem*{conj*}{Conjecture}
\newtheorem*{theo*}{Theorem}
\newtheorem*{ques*}{Question}
\newtheorem*{mi2}{Main Idea}

\newtheorem{ex2}[theo]{Example}

\newtheorem{exer2}[theo]{Exercise}

\newtheorem{alg2}[theo]{Algorithm}

\newcommand{\qq}{{\mathbb{Q}}}                                     
\newcommand{\s}{{\mathbb{S}}}                                      
\newcommand{\U}{{\mathbf{U}}}                                      
\newcommand{\B}{{\mathbf{B}}}                                      
\newcommand{\E}{{\mathbf{E}}}                                      
\newcommand{\dif} {{\operatorname{d}}}                             
\newcommand{\In} {{\,\subseteq\,}}                                 
\newcommand{\id}{{\operatorname{id}}}                              

\newcommand{\comment}[1]{}                                         
\newcommand{\hto}[1]{\overset{#1}{\hookrightarrow}}                

\begin{document}

\maketitle \thispagestyle{empty}


\begin{abstract}
We encode a compact Lie group action on a compact manifold by the Sullivan model of its Borel construction. We then prove that deciding whether this action is almost free is $\mathbf{NP}$-hard.
\end{abstract}


\section*{Introduction}

Computational complexity theory intends to analyse how time or space consuming an optimal algorithm need to be in order to solve a given problem. The theory has vastly thrived throughout the years, with one of its origins certainly lying in the idea of casting problems into different complexity classes. Most prominent amongst the latter are the classes $\mathbf{P}$, which describes all the problems for which there is a polynomial-time solving algorithm, and the class $\mathbf{NP}$, which consists of those problems that may at least be verified in polynomial time. Clearly $\mathbf{P}\In \mathbf{NP}$, however, it is the common belief that several problems in $\mathbf{NP}$ are much harder to solve than the problems in $\mathbf{P}$. Known algorithms typically run at exponential costs.

Moreover, a problem in $\mathbf{NP}$ is said to be \emph{$\mathbf{NP}$-complete} if any
other problem in $\mathbf{NP}$ can be reduced to it in polynomial time. Finally, in
increasing order of difficulty, a problem not necessarily in $\mathbf{NP}$, is \emph{$\mathbf{NP}$-hard}
if, again, any problem in $\mathbf{NP}$ can be reduced to it in polynomial time.

A whole variety of problems stemming from completely different areas of mathematics and computer science have been found to be $\mathbf{NP}$-hard. Just to name a few most prominent ones, we mention the \emph{knapsack problem} and the \emph{subset sum problem}, the \emph{Hamilton circuit problem} and the \emph{travelling salesman problem}, the \emph{satisfiability problem} and the \emph{graph colouring problem}.
(For an introduction to the subject see \cite{Sip06}.)

\vspace{5mm}

Also in the field of algebraic topology problems were considered from the point of view of computational complexity. One area having obtained particular interest is rational homotopy theory. Here, one classical way to transcribe topological spaces to algebraic structures---thus making them accessible to computational complexity theory---is to encode their rational structure by a so-called \emph{Sullivan model}.
Indeed, rational homotopy theory permits a categorical translation from topology to algebra at the expense of losing torsion information, whilst at the same time permitting concrete computations. (As a reference to this theory we recommend the textbook \cite{FHT01}.)

Using this approach several topological problems were shown to be $\mathbf{NP}$-hard. In \cite{Ani89} it is shown that computing the rational homotopy groups $\pi_*(X)\otimes \qq$ of a simply-connected CW-complex $X$ is $\mathbf{NP}$-hard. So is the problem of whether a simply-connected space $X$ with $\dim \pi_*(X)\otimes \qq<\infty$ also has finite-dimensional rational cohomology (cf.~\cite[Theorem 1, p.~90]{LM00}). In the same article it was shown that for formal spaces, i.e.~for spaces for which the rational homotopy type can be formally derived from the rational cohomology algebra, the computation of Betti numbers, of the cup-length and of the rational Lusternik--Schnirelmann category are $\mathbf{NP}$-hard problems. In \cite{GL03} it is shown that the computation of Betti numbers of a simply-connected space with both finite-dimensional rational homotopy and finite-dimensional rational homotopy, a \emph{rationally elliptic space}, is $\mathbf{NP}$-hard.
In \cite{Ama11b} it was proved that even for rationally elliptic spaces the computation of the rational cup-length and the rational Lusternik--Schnirelmann category is $\mathbf{NP}$-hard.

\vspace{5mm}

In this article we consider a problem of a different flavour within this context of computational complexity theory. We consider actions of compact Lie groups on compact manifolds. Recall that such an action is called \emph{free} if all the stabiliser groups $G_x$ (for $x\in M$) are trivial. If these groups are (not necessarily trivial but) still finite, i.e.~$|G_x|<\infty$ for all $x\in M$, then the group action is called \emph{almost free}. Almost free group actions play important roles in several contexts like, for example, the famous toral rank conjecture.

One way to encode the action $G\curvearrowright M$ as a topological space is the so-called \emph{Borel construction} $M_G=(M\times_G \E G)=(M\times \E G)/G$. In particular, this yields the \emph{equivariant cohomology} $H_G^*(M):=H^*(M_G)$ of the $G$-action on $M$. As a transition from topology to algebra we choose the Sullivan model of $M_G$, which results as a twisted product construction of the minimal models of base and fibre in the bundle $M\hto{} M_G\to \B G$.

Using this natural encoding of the action we can prove
\begin{main}\label{theoA}
The decision problem whether the action is almost free or not is $\mathbf{NP}$-hard.
\end{main}

At the end of Section \ref{sec01} we explain what a verifier for this decision problem should be. In this sense we can prove
\begin{maincor}\label{corB}
The decision problem is even $\mathbf{NP}$-complete.
\end{maincor}


\section{Proving the result}\label{sec01}

Let $G$ be a compact Lie group acting on a compact manifold $M$. The main tool which makes deciding almost freeness of this action possible at all is the following famous result (cf.~\cite[Theorem 7.7., p.~276]{FOT08}) coming out of the famous localisation results in equivariant cohomology.
\begin{theo}[Hsiang]\label{theo03}
The action is almost free if and only if $H_G(M;\qq)$ is finite-dimensional.
\end{theo}
As described in the Introduction, we use the Sullivan model of the Borel construction $M_G$ as an encoding of our problem. This model is formed as the model of the fibration
\begin{align*}
M\hto{} M_G \to \B G
\end{align*}
For this we choose a minimal Sullivan model of fibre and base---being the model of a compact Lie group the model of the base is necessarily freely generated in odd degrees---and obtain a twisted differential on the tensor product---the twisting reflecting the complexity of the action (see \cite[Corollary, p.~199]{FHT01}).

Before we can proceed with the actual proof, we need to cite the subsequent lemma from \cite[Lemma 1.2, p.~152]{KZ04}.
\begin{lemma}\label{lemma01}
Let $G$ be a compact Lie group acting differentiably on manifolds $X$ and $Y$. Suppose that the action of $G$ on $X$ is transitive and that the diagonal action of $G$ on $X\times Y$ is free. Then for any $x\in X$ the action of the isotropy group $G_x$ on $Y$ is free and the quotient spaces $(X\times Y)/G$ and $Y/G_x$ are canonically diffeomorphic.
\end{lemma}

\vspace{5mm}

The proof of Theorem \ref{theoA} proceeds as follows. In view of Theorem \ref{theo03} we consider the Sullivan model of the Borel construction and we show that determining whether its cohomology is finite-dimensional is NP-hard.

In \cite{LM00} it is shown that deciding whether the cohomology is finite dimensional is already NP-hard for the following class of finitely-generated Sullivan algebras.
\begin{align}\label{eqn01}
(\Lambda (\langle x_1,\ldots, x_r\rangle\oplus \langle y_1,\ldots, y_s\rangle),\dif)
\end{align}
with $\deg x_i=2$, $\deg y_i$, $\deg y_{(a,b)}=2k-3$, $\dif x_i=0$, $\dif y_{(a,b)}=\sum_{l=1}^k x_a^{k-l} x_b^{l-1}$.

Indeed, in this article the graph colouring problem for undirected, simple, connected, finite graphs is considered. This problem is polynomially reduced to the question whether a Sullivan algebra of the type above has finite dimensional cohomology. For this, vertices are encoded as $x_i$, edges correspond to the $y_i$. An edge between the vertices $a$ and $b$ yields the relation $\dif y_{(a,b)}=\sum_{l=1}^k x_a^{k-l} x_b^{l-1}$.
For $k\geq 3$ it is shown that the algebra is elliptic, i.e.~has finite dimensional cohomology, if and only if the graph is not $k$-colourable.

Thus in order to prove Theorem \ref{theoA} it remains to see that the class of algebras \eqref{eqn01} can be realised by the model of the Borel construction of certain torus actions on compact manifolds. This will show that solving the decision problem in Theorem \ref{theoA} is actually harder than the graph colouring problem---and $\mathbf{NP}$-hard itself.

In order to simplify notation, we shall actually make an index shift $k\mapsto k-1$, i.e.~the space we construct for a fixed $k\geq 2$ will realise $(k+1)$-colourability of the corresponding graph. In other words, we shall realise the algebras
\begin{align}\label{eqn02}
(\Lambda (\langle x_1,\ldots, x_r\rangle\oplus \langle y_1,\ldots, y_s\rangle),\dif)
\end{align}
with $\deg x_i=2$, $\deg y_i$, $\deg y_{(a,b)}=2k-1$, $\dif x_i=0$, $\dif y_{(a,b)}=\sum_{l=0}^k x_a^{k-l} x_b^{l}$.

\begin{proof}[\textsc{Proof of Theorem \ref{theoA}}]
The idea of the proof is to construct a homogeneous space which has the rational type of a product $\s^{2k-1} \times \stackrel{(s)}{\ldots} \times \s^{2k-1}$ together with an action of an $r$-torus $T^r$. Abstractly, this space will be given as
\begin{align*}
B_1\times \ldots \times B_s=B^{\times s}
\end{align*}
together with an action of $T^r$ (which does not respect the product splitting). Each factor $B=B_i$ is a homogeneous space of the form
\begin{align*}
\frac{\U(k)\times \stackrel{(k+2)}{\ldots}\times \U(k)}{\U(k-1)\times \U(k)\times \stackrel{(k+1)}{\ldots}\times \U(k)}
\end{align*}
i.e.~a diffeomorphism sphere $B\cong\s^{2k-1}$.

We construct this homogeneous space together with a torus action. For this we consider
\begin{align*}
\U(k)\times \stackrel{(k+2)}{\ldots}\times \U(k)
\end{align*}
together with the action of
\begin{align*}
(T^k)_1\times \stackrel{(k+2)}{\ldots} \times (T^k)_{k+2}\times \U(k-1)\times (\U(k))_2 \times \stackrel{(k)}{\ldots}\times(\U(k))_{k+2}
\end{align*}
defined in the following way. The $(T^k)^{k+2}$ torus acts from the left in the standard way. Each $T^k$ is the standard maximal torus of the respective $\U(k)$. The remaining factors act from the right. The inclusions are as follows: The group $\U(k-1)$ is given by diagonal inclusion into  $\U(k)\times \stackrel{(k+2)}{\ldots}\times \U(k)$ where the inclusion on each factor is the standard one. The group $\U(k)_i$ is included in the standard way into the first and the $i$-th factor of the product $\U(k)\times \stackrel{(k+2)}{\ldots}\times \U(k)$; i.e.~an element is mapped as $u\mapsto (u)_1 \times \id \times \ldots \times \id \times (u)_i \times \id \times \ldots \times \id$.

We then make two claims:
\begin{enumerate}
\item Via these inclusions, $\U(k-1)\times (\U(k))_2 \times \stackrel{(k)}{\ldots}\times(\U(k))_{k+2}$ is a subgroup of $\U(k)\times \stackrel{(k+2)}{\ldots}\times \U(k)$, and each $B_i$ is diffeomorphic to $\s^{2k-1}$. The volume form corresponds to $x_1-x_2- \ldots - x_{k+2}$, where the $x_i$ are the volume forms of the $\U(k)$-nominator-factors.
\item There is a torus action on $B_1\times \ldots \times B_s$ such that the Borel construction of the action is as required.
\end{enumerate}

\textbf{ad (1).} It is easy to see that for $k\geq 1$ the inclusions of the factors yield an injective group homomorphism on the direct product; and we have a compact homogeneous space. Successively applying Lemma \ref{lemma01} to the $(\U(k))_i$ we see that
\begin{align*}
&{\U(k)\times \stackrel{(k+2)}{\ldots}\times \U(k)}/{\U(k-1)\times \U(k)\times \stackrel{(k+1)}{\ldots}\times \U(k)}\\\cong & \U(k)/\U(k-1)\\\cong& \s^{2k-1}
\end{align*}

Forming the model of this homogeneous space
\begin{align*}
(\Lambda (V_{\B \U(k-1)}\oplus V_{\B \U(k)} \oplus \ldots \oplus V_{\B \U(k)}) \otimes \Lambda (V_{\U(k)}\oplus \ldots \oplus V_{\U(k)}),\dif)
\end{align*}
we derive the assertion on the volume form by computing that $\dif x_1=\sum_{1\leq i\leq k+1} y_i$, $\dif x_i=y_i$ for $i\geq 2$ (where $y_i$ for $1\leq i\leq k+1$ is the volume form of a denominator $\U(k)$-factor). Indeed, this yields that $\ker \dif =\langle x_1- x_2- \ldots - x_{k+1}\rangle$.

\textbf{ad (2).} First, we compute the Sullivan model of the Borel construction belonging to the action of $T^k$ on $B_i$. We then embed a $2$-torus into $T^k$ such that the relation $\dif y_{(a,b)}=\sum_{l=0}^k x_a^{k-l} x_b^{l}$ is realised. We then extend this inclusion to the inclusion of $T^r$ into $(T^k)^{k+2}$ such that the algebra \eqref{eqn02} is realised.

The Borel construction $(G/H)_T=\E T\times_T G/H $ fits in the fibration
\begin{align*}
G/H \hto{} \E T\times_T G/H \to \B T
\end{align*}
We now consider the product action of $T$ on $(\E T\times_T G/H) \times T$ (with the standard action on the last factor). This action is free and fits into a fibration which we may compare to the first one. We obtain the commutative diagram
\begin{align*}
\xymatrix{
G\times T/H\ar@{^{(}->}[d] \ar[r] &G/H \ar@{^{(}->}[d]\\
T\backslash G\times T/H\simeq_\qq T \backslash \E T \times G\times T/H  \ar[r]^<<<<<p\ar[d] & T \backslash \E T \times G/H
\ar[d]\\
\B T\ar@{=}[r] & \B T }
\end{align*}
with the obvious projection $p$.

We form the models $(\Lambda V_{\B T} \otimes \Lambda V_{G} \otimes \Lambda V_{\B H},\dif)$ of $(G/H)_T$ and $(\Lambda V_{\B T} \otimes \Lambda V_{G} \otimes \Lambda V_{\B H},\dif)$ of $(G\times T/H)_T$. The morphism of Sullivan algebras induced by the diagram morphism is
\begin{align*}
(\Lambda V_{\B T} \otimes \Lambda V_{G} \otimes \Lambda V_{\B H},\dif) \to (\Lambda V_{\B T} \otimes \Lambda V_{G} \otimes \Lambda V_{\B H} \otimes \Lambda V_T,\dif)
\end{align*}
given by the identities on the respective factors. This shows that a model for the Borel construction $(G/H)_T$ is given by.
\begin{align*}
(\Lambda V_{\B T} \otimes \Lambda V_{G} \otimes \Lambda V_{\B H},\dif)
\end{align*}
with $\dif v=H(\B G\to \B H)v^{+1}-H(\B G \to \B T)v^{+1}$ (as in the biquotient case and under the usual identifications).

\vspace{5mm}

Let us come back to the concrete case of the spheres $B_i$ from above. We specifiy the action of an $r$-torus $T^r$ on $B_1 \times \ldots \times B_s$ as follows. Whenever the Sullivan algebra \eqref{eqn02} we want to encode has a relation of the form $\dif y_{(a,b)}=\sum_{l=0}^k x_a^{k-l} x_b^{l}$ we include the standard rank 2 subtorus $(\s^1)_a\times (\s^1)_b\In T^r$ into the maximal torus $T$ of $B_{(a,b)}$ acting from the left on $B_{(a,b)}$ in the way we shall specify below. The product of all these inclusions then defines the $T^r$-action on $B_1\times \ldots \times B_s$. The model of the Borel construction of this action is the Sullivan algebra \eqref{eqn02} we started with.

It remains to specify the inclusion of such a $2$-torus into one maximal torus of some $B_{(a,b)}$. On each $B_{(a,b)}$ we considered the left action of a $(T^k)^{k+2}=(T^k)_{-1} \times \ldots \times (T^k)_{k+1}$. We describe the inclusion of $\s^1_{a}\times \s^1_{b}$ into this torus via the inclusions into the respective $(T^k)_i$ given as
\begin{align*}
(t_a, \stackrel{(i)}{\ldots}, t_a,t_b,\stackrel{(k-i)}{\ldots}, t_b)
\end{align*}
for $0\leq i\leq k$.

In other words, we never include into the first $T^k$-factor, and we use the standard inclusion of $\s_a$ on $i$ factors, the standard inclusion of $\s_b$ on $k-i$ factors. The first $t_b$ occurs at position $i+1$.

Using Claim (1) and the form of the volume form we derived there, it is now immediate that in the Sullivan model the differential of the volume form of $B_{(a,b)}$ has $-t_a^it_b^{k-i}$ as a summand. Combining all the inclusions, we derive the result.
\end{proof}
\begin{rem}
We remark that the model of the Borel construction (of a torus action) may always be realised by a \emph{free} action once it has finite dimensional cohomology using the classifying space of torus bundles (cf.~\cite[Proposition 7.17, p.~280]{FOT08}).
\end{rem}

\vspace{5mm}

Let us finally deal with the assertion of Corollary \ref{corB}. In \cite{LM00} the finite-dimensionality of the algebra \eqref{eqn01} is verified by the existence of a non-trivial morphism of differential graded algebras to the free polynomial algebra $\qq[z]$ with $\deg z=2$---see \cite[Proof of Theorem 3, p.~90]{LM00}. It is shown that the corresponding graph is $k$-colourable if and only if the algebra is not elliptic if and only if such a morphism (and with it a non-nilpotent element in cohomology) exists. In other words, a proposed colouring of the graph transcribes to prescribing such a morphism on the generators in degree $2$. Given such a morphism on the degree $2$ generators, it is certainly possible  to check  in polynomial time whether it commutes with differentials and hence defines a morphism of Sullivan algebras. This makes this problem from \cite{LM00} $\mathbf{NP}$-complete.

We do exactly the same in our case and verify almost freeness by the existence of a non-trivial morphism to $\qq[z]$ commuting with differentials as above. This then proves Corollary \ref{corB}.



\pagebreak \

\vfill

\begin{center}
\noindent
\begin{minipage}{\linewidth}
\small \noindent \textsc
{Manuel Amann} \\
\textsc{Fakult\"at f\"ur Mathematik}\\
\textsc{Institut f\"ur Algebra und Geometrie}\\
\textsc{Karlsruher Institut f\"ur Technologie}\\
\textsc{Englerstrasse 2}\\
\textsc{76131 Karlsruhe}\\
\textsc{Germany}\\
[1ex]
\textsf{manuel.amann@kit.edu}\\
\textsf{http://topology.math.kit.edu/$21\_54$.php}
\end{minipage}
\end{center}

\end{document}